\documentclass[12pt]{article}
\usepackage{latexsym}
\usepackage{amsmath, amsthm}
\usepackage[all,dvips,ps]{xy}
\usepackage{amsfonts}
\usepackage{mathrsfs}
\usepackage{enumerate}

\newtheorem{prop}{Proposition}[section]
\newtheorem{lem}{Lemma}[section]
\newtheorem{thm}{Theorem}[section]
\newtheorem{cor}{Corollary}[section]

\newtheorem{exa}{Example}[section]

\newcommand{\tra}{\mathrm{trace}\,}
\newcommand{\diag}{\mathrm{diag}\,}
\newcommand{\rank}{\mathrm{rank}\,}

\newcommand{\Rs}{\mathbb{R}}

\newcommand{\thl}{\underline{\theta}}
\newcommand{\thu}{\overline{\theta}}
\newcommand{\thc}{\theta_c}
\newcommand{\Tl}{T^{\leq}_{kl}}
\newcommand{\Te}{T^{=}_{kl}}

\newcommand{\zz}{\tilde{z}}
\newcommand{\ZZ}{\tilde{Z}}

\newcommand{\DD}{\tilde{D}}
\newcommand{\ww}{\tilde{w}}
\newcommand{\ee}{\tilde{e}}
\newcommand{\PP}{\tilde{P}}

\newcommand{\N}{\mbox{null}}
\newcommand{\bz}{{\bf 0} }

\newcommand{\bpr}{{\bf Proof.} \hspace{1 em}}
\newcommand{\epr}{ \\ \hspace*{4.5in} $\Box$ }
\newcommand{\beq}{ \begin{equation} }
\newcommand{\eeq}{ \end{equation} }
\newcommand{\bt}{ \begin{tabular} }
\newcommand{\et}{ \end{tabular} }

\begin{document}

\bibliographystyle{plain}
\title{On Unit Spherical Euclidean Distance Matrices Which Differ in One Entry}
\vspace{0.3in}
        \author{ A. Y. Alfakih
  \thanks{E-mail: alfakih@uwindsor.ca}
  \\
          Department of Mathematics and Statistics \\
          University of Windsor \\
          Windsor, Ontario N9B 3P4 \\
          Canada
}

\date{March 15, 2019. Revised \today}
\maketitle


\vspace{0.1in}

\begin{abstract}
A unit spherical Euclidean distance matrix (EDM) $D$ is a matrix whose entries can be realized
as the interpoint (squared) Euclidean distances of $n$ points on a unit sphere.
In this paper, given such a $D$ and $1 \leq k < l \leq n$, we present a characterization of the set of all unit spherical EDMs
whose entries agree with those of $D$ except possibly with the entry in the $kl$th and $lk$th positions.
As a result, we show  that
this set can be discrete, consisting of  one or two elements, or it can be continuous. The results are derived using two
alternative approaches, the second of which is based on Cayley-Menger matrices.
\end{abstract}

\section{Introduction}
An $n \times n$ matrix $D=(d_{ij})$ is said to be a Euclidean distance matrix (EDM) if there exist points $p^1,\ldots,p^n$ in some Euclidean
space such that
\[
d_{ij}= || p^i - p^j ||^2 \mbox{ for all } i,j=1,\ldots, n,
\]
where $||x||$ denotes the Euclidean norm of $x$.
$p^1,\ldots,p^n$ are called the {\em generating points} of $D$ and the dimension of their affine span is
called the {\em embedding dimension} of $D$.  If the  generating points of an EDM $D$ lie on a
sphere of radius $\rho$, then $D$ is called {\em spherical} and $\rho$ is referred to as the radius of $D$.
A spherical EDM of unit radius is called a {\em unit spherical} EDM.

Let $E^{kl}$ denote the $n \times n$ symmetric matrix with 1's in the $kl$th and $lk$th positions
and zeros elsewhere. Let $D$ be an EDM and let $l_{kl} \leq 0$ and $u_{kl} \geq 0$ be the two scalars  such that
$D + t E^{kl}$ is an EDM if and only if $l_{kl} \leq t \leq u_{kl}$. That is, $D$ remains an EDM if its
entry in the $kl$th and $lk$th positions varies between $d_{kl}+l_{kl}$ and $d_{kl} + u_{kl}$,
while keeping all other entries unchanged. The entry $d_{kl}$ is said to be {\em unyielding} if $u_{kl} = l_{kl} = 0$
and it is said to be {\em yielding} if $u_{kl} \neq l_{kl}$. The recent paper \cite{alf18} presents a characterization of
the yielding entries of an EDM $D$ and derives simple closed-form expressions of  $u_{kl}$ and $l_{kl}$
for each yielding entry $d_{kl}$.

In this paper we extend the work in \cite{alf18} and focus on unit spherical EDMs.
Let $D$ be a given unit spherical EDM  and let $d_{kl}$ be a yielding entry of $D$.
Let
\beq \label{defTless}
\Tl = \{ t \in [l_{kl},u_{kl}] : D + t E^{kl} \mbox{ is a spherical EDM of radius $\rho \leq 1$} \}.
\eeq
Obviously, $0 \in \Tl$.
 We characterize those yielding entries $d_{kl}$ for which
$\Tl \neq \{0\}$ in terms of Gale transform of the
generating points of $D$ and in terms of vector $w$ (defined in (\ref{defw}) below).
Moreover, for such entries we derive closed-form expressions of $\Tl$.
It is worth pointing out that if $t \in [l_{kl},u_{kl}] \backslash \Tl$, then $D +t E^{kl}$ may be spherical, however,
its radius is $> 1$.

Now let
\beq \label{defTless}
\Te = \{ t \in \Tl : D + t E^{kl} \mbox{ is a unit spherical EDM} \}.
\eeq

Unlike $\Tl$, which is a closed subinterval of $[l_{kl},u_{kl}]$, $\Te$ is not necessarily convex.
The main result of this paper is a closed-form expression of $\Te$. As a result, we show that
$\Te$ can have one of three possible forms: First,
it can be a subinterval of $\Tl$, i.e., $D$ can have a continuum of unit spherical EDMs which differs from it only
in the entry in the $kl$th and $lk$th positions. Second, it can consist of two points, one of which is obviously 0,
i.e., $D$ can have exactly one other unit spherical EDM which differs from it only
in the entry in the $kl$th and $lk$th positions. Third, it can be a singleton, i.e., $\Te=\{0\}$, in which case
there does not exist a unit spherical EDM which differs from $D$ only
in the entry in the $kl$th and $lk$th positions.
Our results are derived using two alternative approaches, the second of which is based on Cayley-Menger matrices.
As a by-product of our characterizations of
$\Tl$ and $\Te$ we obtain some other related results.

The remainder of this paper is organized as follows. In Section 2 we review the mathematical background needed in later
sections.  The characterizations of $\Tl$ and $\Te$ are presented in Sections 3 and 4 respectively. Finally, in Section 5,
we rederive the results of Sections 3 and 4 using an alternative approach based on Cayley-Menger matrices.

\subsection{Notation}
We collect here the notation used throughout  the paper. $E^{kl}$ denotes the $n \times n$ matrix with 1's
 in the $kl$th and $lk$th positions and 0's elsewhere. $e^i$ denotes the standard unit vector in $\Rs^n$ and thus
 $E^{kl} = e^k (e^l)^T +  e^l ( e^k)^T$.
 $e$ and $\ee$ denote, respectively, the vectors of all 1's in $\Rs^n$ and $\Rs^{n+1}$. $E = e e^T$
 denotes the $n \times n$ matrix of all 1's and $I_n$ denotes
the identity matrix of order $n$. The zero matrix or zero vector of appropriate dimension is denoted  by $\bz$. For a symmetric
matrix $A$, we mean by $A \succeq \bz$ that $A$ is positive semidefinite.
$A^{\dag}$ denotes the Moore-Penrose inverse of $A$ and
$\N(A)$ denotes the null space of $A$.
Finally, $\backslash$ denotes the set theoretic difference.

\section{Preliminaries}

In this section, we present some known and some new results concerning EDMs, spherical EDMs, Gale matrices and
Gram matrices that will be needed in the sequel. Also, we present some of the main results of \cite{alf18}.
For a comprehensive treatment of EDMs see the monograph \cite{alfm18}.

\subsection{EDMs and Gale Matrices}

Let $e$ denote the vector of all 1's in $\Rs^n$.
For a matrix $A$, we use $A \succeq \bz$ to indicate that $A$ is symmetric positive semidefinite.
Let $D$ be an $n \times n$ real symmetric matrix whose diagonal entries are all 0's. Then
it is well known \cite{sch35,yh38, gow85,cri88}
that $D$ is an EDM if and only if $D$ is negative semidefinite on $e^{\perp}$, the orthogonal complement of $e$ in $\Rs^n$.
In other words, if $s$ is a vector
in $\Rs^n$ such that $e^Ts= 1$, then $D$ is an EDM if and only if
\beq  \label{defB}
B= -\frac{1}{2} (I-e s^T ) D (I - s e^T) \succeq \bz,
\eeq
in which case, the embedding dimension of $D$ is given by rank$(B)$.

Let $D$ be an $n \times n$ EDM of embedding dimension $r$ and let $B$ be as defined in (\ref{defB}). Then $B$ can be factorized as  $B =PP^T$
where $P$ is $n \times r$. Consequently, $p^1, \ldots, p^n$, the generating points of $D$, are given by the rows of $P$.
That is,
\beq
P = \left[ \begin{array}{c} (p^1)^T \\ \vdots \\ (p^n)^T \end{array} \right].
\eeq
As a result, $P$ is called a {\em configuration matrix} of $D$.
The following three observations are immediate. First,
$P$ has full column rank, i.e., rank$(P) = r$. Second,
$B$ is the Gram matrix of the generating points of $D$ (or the Gram matrix of $D$ for short).
Third, $P^Ts=\bz$ since $Bs=\bz$.

It is well known \cite{gow85} that if $D$ is a nonzero EDM, then $e$ lies in the column space of $D$. Hence, there exists $w$ such that
\beq  \label{defw}
D w = e.
\eeq
We assume that $w=D^{\dag}e$, i.e., $w$ is perpendicular to $\N(D)$. Vector
$w$ plays a key role in this paper.

 Different choices of vector $s$ in (\ref{defB}) amount to different choices of the origin.
Two choices of $s$ are of particular interest in this paper:

First,  $s=e/n$. In this case, let
\beq \label{defJ}
  J:= I - ee^T/n.
\eeq
Hence, the corresponding Gram matrix is given by $B=-JDJ/2$, and consequently the origin is fixed  at the centroid of
the generating points of $D$ since $Be=\bz$.

Second, $s=2w$, where $Dw=e$. Assume that $D$ is a unit spherical EDM. Then, by Theorem \ref{thmsEDM} below,
 $2 e^T w = 1$. Consequently, in this case, the
corresponding Gram matrix is given by $B'= E - D/2$ and satisfies $B'w=\bz$.

Assume that $r \leq n-2$ and let
\beq
Z = \left[ \begin{array}{c} (z^1)^T \\ \vdots \\ (z^n)^T \end{array} \right]
\eeq
 be the $n \times (n-r-1)$ matrix whose columns form a basis of
\beq
\N(\left[ \begin{array}{c} P^T \\ e^T \end{array} \right])= \N(\left[ \begin{array}{c} B \\ e^T \end{array} \right]).
\eeq
$Z$ is called a {\em Gale matrix} of $D$ and $z^1, \ldots, z^n \in \Rs^{n-r-1}$ are called
{\em Gale transforms} of $p^1, \ldots, p^n$. The notion of
Gale transform \cite{gal56, gru67} is well known and widely used in the theory of polytopes.
The components $w_k$ and $w_l$ of  $w$ together with the Gale transforms $z^k$ and $z^l$ play the crucial role
in the characterizations of the sets $\Tl$ and $\Te$.

\subsection{Spherical EDMs}

The following known characterizations of spherical EDMs
 \cite{aw02,gow82,gow85,thw96,neu81} are needed in the sequel.

\begin{thm} \label{thmsEDM}
Let $D$ be an $n \times n$ EDM of embedding dimension $r$ and let $Dw =e$. Let $Z$ and $P$ be, respectively, a
Gale matrix and a configuration matrix of $D$ such that $P^T e =\bz$.
If $r = n-1$, then $D$ is  spherical. Otherwise, if $r \leq n-2$, then
the following statements are equivalent:
\begin{enumerate}
\item $D$ is spherical,
\item the columns of $Z$ form a basis of $\N(D)$
\item $\rank(D)=r+1$.
\item there exists $a \in \Rs^r$ such that ${\displaystyle P a = \frac{1}{2} J \diag(PP^T)}$,
in which case, the generating points of $D$ lie on a sphere centered at $a$ and has radius
\[
 \rho = \left(a^Ta + \frac{e^T D e}{2 n^2} \right)^{1/2}.
\]
\item $e^Tw >0$, in which case, the radius of $D$ is given by  ${\displaystyle \left(\frac{1}{2 e^T w } \right)^{1/2}}$.
\item there exists a scalar $\beta$ such that  $\beta e e^T - D \succeq \bz$,
 and the minimum value of such a $\beta$ is $2 \rho^2$, where $\rho$ is the radius of $D$.
\end{enumerate}
\end{thm}

An interesting class of spherical EDMs is that of regular EDMs. A spherical EDM $D$ is {\em regular} if the center of
the sphere containing the generating points of $D$ coincides with their centroid. That is, if $a=\bz$ in part 4 of Theorem \ref{thmsEDM}
and thus $\rho^2= e^TDe/(2n^2)$ in this case. It is known \cite{ht93} that an EDM $D$ is regular if and only if its
Perron eigenvector is $e$. In other words, iff $De=(e^TDe/n) e$. Consequently, $w = e /(2n \rho^2)$
for regular EDMs.

 Since EDMs are either spherical or nonspherical, many characterizations of spherical EDMs provide at the same
 time characterizations of nonspherical EDMs.
The most useful characterizations of nonspherical EDMs for our purposes are given in the following theorem.

\begin{thm}[Gower \cite{gow82,gow85}]  \label{thmnsEDM}
Let $D$ be an $n \times n$ EDM of embedding dimension $r$ and let $Dw =e$.
Then the following statements are equivalent:
\begin{enumerate}
\item $D$ is nonspherical,
\item  $e^Tw=0$.
\item $\rank(D)=r+2$.
\end{enumerate}
\end{thm}

\subsection{The Gram Matrices $B$ and $B'$}

Let $B$ and $B'$ be the two Gram matrices corresponding to a unit spherical EDM $D$ such that $Be=\bz$ and $B'w=\bz$. In other words,
let $B=-JDJ/2$ and $B'=E-D/2$ and let $B$ and $B'$ be factorized as $B=PP^T$ and $B'=P'P'^T$.
Define
\beq  \label{defSS'}
S= P(P^TP)^{-1}    \mbox{ and } S'=P'(P'^TP')^{-1},
\eeq
and let $(s^i)^T$ and $(s'^i)^T$ denote the $i$th rows of $S$ and $S'$ respectively. The vectors ${s^i}$ and ${s'^i}$ play a key role
in this paper.

It is easy to see that $SS^T = P(P^TP)^{-2} P^T = B^{\dag}$, where $B^{\dag}$
is the Moore-Penrose inverse of $B$. Similarly, $S'S'^T = B'^{\dag}$.
In this subsection we derive some useful properties of $B^{\dag}$ and $B'^{\dag}$.

The following two theorems are slight generalizations of a theorem of Styan and Subak-Sharpe \cite{ss97}.

\begin{thm} \label{thmssBw}
Let $D$ be an $n \times n$ unit spherical EDM and let $B'=E-D/2$. Then
\[
B'^{\dag} = -2 D^{\dag} + \frac{2}{w^Tw} (D^{\dag} w w^T + ww^T D^{\dag} - \frac{w^T D^{\dag}w}{w^Tw} ww^T),
\]
where $w=D^{\dag} e$.
\end{thm}

\bpr
Let $r$ be the embedding dimension of $D$. Let  $\ZZ = w$ if $r=n-1$ and $\ZZ=[w \;\; Z]$ if $r \leq n-2$, where
$Z$ is a Gale matrix of $D$. Then the columns of  $\ZZ$ form a basis of $\N(B')$ and thus
$B'^{\dag} B' = I - \ZZ (\ZZ^T \ZZ)^{-1} \ZZ^T$. Therefore,

\beq  \label{eqstyan}
B'^{\dag} D = 2 B'^{\dag} e e^T - 2 I + 2 \ZZ (\ZZ^T \ZZ)^{-1} \ZZ^T.
\eeq
Assume that $r \leq n-2$. The case where $r=n-1$ is similar and easier.
Hence, by multiplying (\ref{eqstyan}) from the right by $D^{\dag}$ and by using the facts that
$D D^{\dag} = I - Z (Z^TZ)^{-1} Z^T$ and $B'^{\dag} Z = \bz$ , we have
\beq  \label{eqstyan2}
B'^{\dag} = 2 B^{\dag} e w^T - 2 D^{\dag} + 2 \ZZ (\ZZ^T \ZZ)^{-1} \ZZ^T D^{\dag},
\eeq

The fact $B'^{\dag}w = 0$ implies that
\[
  B'^{\dag} e = \frac{1}{w^Tw} (D^{\dag} w   -\ZZ (\ZZ^T \ZZ)^{-1} \ZZ^T D^{\dag}w ).
\]
Substituting $B'^{\dag} e$ into (\ref{eqstyan2}) we have
\[
B'^{\dag} = -2D^{\dag} + \frac{2}{w^Tw} (D^{\dag} w w^T  -\ZZ (\ZZ^T \ZZ)^{-1} \ZZ^T D^{\dag}w w^T)  + 2 \ZZ (\ZZ^T \ZZ)^{-1} \ZZ^T D^{\dag}.
\]
Now $Z^T w = Z^T D^{\dag}e = \bz$. Thus, $(\ZZ^T \ZZ)^{-1} = \left[ \begin{array}{cc} (w^Tw)^{-1} & \bz  \\ \bz &  (Z^TZ)^{-1} \end{array} \right]$.
Therefore,
$\ZZ (\ZZ^T \ZZ)^{-1} \ZZ^T D^{\dag}= (w^Tw)^{-1} ww^T D^{\dag}$ and thus the result follows.
\epr

The argument used in the proof of the previous theorem can also be used to prove the following theorem.

\begin{thm}[\cite{bb07}]  \label{thmssB}
Let $D$ be a unit spherical EDM  and let $B=-JDJ/2$. Then
\[
B^{\dag} = -2 D^{\dag} + 4 ww^T,
\]
where $w= D^{\dag} e$.
\end{thm}

The following corollaries, which will be used in the sequel, are immediate consequences of Theorems \ref{thmssBw} and \ref{thmssB}.

\begin{cor} \label{corBiBi'w0}
Let $D$ be a unit spherical EDM and  let $w= D^{\dag} e$.
Let $s^i$ and $s'^i$ be as defined in (\ref{defSS'}) and assume that $w_k=w_l=0$. Then
\[
(s^k)^T s^l \pm ||s^k|| \; ||s^l|| = (s'^k)^T s'^l \pm ||s'^k|| \; ||s'^l||.
\]
\end{cor}

\bpr
This is immediate since $B^{\dag}_{kk} = B'^{\dag}_{kk}= -2 D^{\dag}_{kk}$, $B^{\dag}_{ll} = B'^{\dag}_{ll}= -2 D^{\dag}_{ll}$ and
$B^{\dag}_{kl} = B'^{\dag}_{kl}= -2 D^{\dag}_{kl}$.

\begin{cor} \label{corBiBi'c}
Let $D$ be a unit spherical EDM and  let $w= D^{\dag} e$.
Let $s^i$ and $s'^i$ be as defined in (\ref{defSS'}) and assume that $w_k=cw_l$, where $w_k \neq 0$. Then
\[
D^{\dag}_{kk} + c^2 D^{\dag}_{ll}-2 c D^{\dag}_{kl} = - \frac{1}{2} ||s^k - c s^l||^2 = - \frac{1}{2} ||s'^k - c s'^l||^2.
\]
\end{cor}

\bpr
Let $x = e^k - c e^l$, where $e^i$ is the $i$th standard unit vector in $\Rs^n$.
Then $w^T x = 0$. The result follows from Theorems \ref{thmssBw} and \ref{thmssB} since
$x^T B^{\dag}x = -2 x^T D^{\dag} x = x^T B'^{\dag}x $.
\epr

We end this subsection with the following lemma which is a simple corollary of Theorem \ref{thmssB}.

\begin{lem} \label{lemDdagBdag}
Let $D$ be a unit spherical EDM and let $B=-JDJ/2$. Then
\[
\DD^{\dag} =\left[ \begin{array}{cc} 0 & e^T \\ e & D \end{array} \right]^{\dag}
 = \left[ \begin{array}{cc} -2 & 2 w^T \\ 2 w & -  B^{\dag}/2 \end{array} \right],
\]
where $w= D^{\dag} e$.
\end{lem}

\subsection{Yielding Entries of an EDM}

In this subsection we review some of the main results which we need from \cite{alf18}.
Vectors $u$ and $v$ in $\Rs^n$ are {\em parallel} if $u = c v$ for some nonzero scalar $c$. Thus,
if $u=v=\bz$, then $u$ and $v$ are parallel. The following proposition characterizes the eigenvalues of
rank-two symmetric matrices.

\begin{prop}[\cite{alf18}] \label{prop1}
Let $\Psi=a b^T + b a^T$, where $a$ and $b$ are two nonzero, nonparallel vectors in $\Rs^r$, $r \geq 2$. Then $\Psi$ has
exactly one positive eigenvalue $\lambda_1$ and one negative eigenvalue $\lambda_r$, where
\[
\lambda_1 = a^Tb + ||a|| \; ||b||   \mbox{ and } \lambda_r = a^Tb - ||a|| \; ||b||.
\]
\end{prop}

The yielding entries of an EDM $D$ are characterized in the following theorem.

\begin{thm}[\cite{alf18}] \label{thmfrom2}
Let $D$ be an $n \times n$ EDM of embedding dimension $r$. If $r=n-1$, then each entry of $D$ is yielding.
On the other hand, if $r \leq n-2$, let $z^1,\ldots, z^n$ be Gale transforms of the generating points of $D$.
Then entry $d_{kl}$ is yielding if and only if $z^k$ is  parallel to $z^l$; i.e., iff
there exists a nonzero scalar $c$ such that $z^k = c z^l$.
\end{thm}

Let $d_{kl}$ be a given yielding entry of $D$.
Before presenting a characterization of the yielding intervals of $d_{kl}$,
it is convenient to define the following quantities which will be used throughout the paper.
Let $B=-JDJ/2$ and let  $s^k$ and $s^l$ be as defined in (\ref{defSS'}). Define

\beq \label{defthul}
 \thl= \frac{2}{(s^k)^T s^l - ||s^k|| \; ||s^l||} \;\; \mbox{ and } \;\; \thu= \frac{2}{(s^k)^T s^l + ||s^k|| \; ||s^l||}.
\eeq
and
\beq  \label{defthc}
 \thc = \frac{-4c}{||s^k-cs^l||^2}.
\eeq

\begin{thm}[\cite{alf18}] \label{thmfrom2int}
Let $D$ be an $n \times n$ EDM of embedding dimension $r$ and let $B=-JDJ/2$.
Assume that the entry $d_{kl}$ is yielding.
If $r=n-1$ or if $r \leq n-2$ and $z^k=z^l=\bz$, then the yielding interval of $d_{kl}$ is given by
\[
[l_{kl},u_{kl}] = \left[  \thl \; , \; \thu  \right],
\]
where $\thl$ and $\thu$ are defined in (\ref{defthul}).

On the other hand, if $r \leq n-2$ and $z^k=c z^l \neq \bz$ where $c$ is a nonzero scalar, then
the yielding interval of $d_{kl}$ is given by
\[
[l_{kl},u_{kl}] = \left\{ \begin{array}{ll} \left[  \thc \; , \; 0  \right] & \mbox{ if } c > 0, \\
                                        \left[ 0 \; , \; \thc  \right]  & \mbox{ if } c < 0, \end{array} \right.
\]
where $\thc$ is defined in (\ref{defthc}).
\end{thm}

\section{Characterizing $\Tl$}

Let $D$ be an $n \times n$  unit spherical EDM of embedding dimension $r$.
Assume that $d_{kl}$ is a yielding entry of $D$ with yielding interval $[l_{kl}, u_{kl}]$.
Recall that
\[
\Tl=\{t \in [l_{kl},u_{kl}]  : D+tE^{kl} \mbox{ is a spherical EDM of radius } \rho \leq 1 \}.
\]
As we show in this section, $\Tl$ can be a singleton, i.e., $\Tl=\{0\}$, or it can be an interval of nonzero length.
In the latter case, $\Tl$ may be equal to, or may be a proper subinterval of,
 the yielding interval $[l_{kl}, u_{kl}]$.

Let $B'=E-D/2 = P'P'^T$.
Throughout this paper, it is convenient to define \footnote{The interpretation of $\ZZ$  in given in Section 5.}

\beq \label{defZZ}
\ZZ= \left\{ \begin{array}{cl} w  & \mbox{ if } r= n-1, \\
                                        \left[w \; Z  \right]  & \mbox{ if } r \leq n-2, \end{array} \right.
\eeq
where $w=D^{\dag}e$ and $Z$ is a Gale matrix of $D$. As a result, $\N(P'^T) \cap \N(\ZZ^T) = \{0\}$.

\begin{prop} \label{prop2}
Let $D$ be an $n \times n$ spherical EDM  of embedding dimension $r$,  and let $P'$ be a configuration matrix
of $D$ such that $P'^Tw=\bz$. Further, let $\ZZ$ be as defined in (\ref{defZZ}) and
let $(\zz^i)^T$ denote the $i$th row of $\ZZ$.
\begin{enumerate}
\item If $\zz^l= \bz$, then $p'^l \neq \bz$.
\item If $\zz^k = c \zz^l$ for some nonzero scalar $c$, then $p'^k-cp'^l \neq \bz$.
\end{enumerate}
\end{prop}

\bpr
To prove part(1), assume to the contrary that $\zz^l = \bz$ and $p'^l=\bz$.
Recall that $e^i$ denotes the $i$th standard unit vector in $\Rs^n$.
Then $e^l \in \N(P'^T) \cap \N(\ZZ^T)$, a contradiction.
Similarly, to prove part (2), assume to the contrary that $p'^k-cp'^l = \bz$ and
let $x$ be the vector in $\Rs^n$ with $1$ and $(-c)$ in the $k$th and $l$th positions respectively and 0's elsewhere.
Then $x \in \N(P^T) \cap \N(\ZZ^T)$, a contradiction.
\epr

An immediate consequence of Proposition \ref{prop2} is that
if $\zz^k=\zz^l=\bz$, then $p'^k$ is not parallel to $p'^l$, i.e., there does not exist a nonzero scalar $c$
such that $p'^k=cp'^l$.

\begin{lem}  \label{lembasic}
Let $B'=P'P'^T$. Then
$t \in \Tl$ if and only if
\[
 \left[ \begin{array}{rr} 2(P'^TP')^2 - t \;( p'^k (p'^l)^T + p'^l (p'^k)^T) &
                         -t \; (p'^k (\zz^l)^T + p'^l (\zz^k)^T) \\ -t \; (\zz^k (p'^l)^T + \zz^l (p'^k)^T) &
                         - t \; (\zz^k (\zz^l)^T + \zz^l (\zz^k)^T) \end{array} \right] \succeq \bz,
\]
where $(\zz^i)^T$ is the $i$th row of $\ZZ$ defined in (\ref{defZZ}).
\end{lem}

\bpr
Evidently,
the $n \times n$ matrix $Q= [P' \; \ZZ ]$ is nonsingular.
By Thereom \ref{thmsEDM}, $D + t E^{kl}$ is a spherical
EDM of radius $\rho \leq 1$ if and only if
$2E - (D + t E^{kl}) \succeq \bz$ iff
$Q^T (2B' - t E^{kl}) Q $ and the result follows
\epr

We first characterize the case where $\Tl$ is a singleton, i.e., $\Tl=\{0\}$.

\begin{thm} \label{thmTl0}
Let $D$ be an $n \times n$ unit spherical EDM and let $d_{kl}$ be a yielding entry of $D$.
Let $(\zz^i)^T$ denote the $i$th row of $\ZZ$ defined in (\ref{defZZ}).
Then $\Tl=\{0\}$ if and only if $\zz^k$ is not parallel to $\zz^l$; i.e., iff
there does not exist a nonzero scalar $c$ such that $w_k = c w_l$ and $z^k = c z^l$.
\end{thm}

\bpr
Assume that $\zz^k = c \zz^l$ for some nonzero scalar $c$.
Then $\zz^k (\zz^l)^T + \zz^l (\zz^k)^T = 2 c \zz^l (\zz^l)^T$ and
$p'^k (\zz^l)^T + p'^l (\zz^k)^T = (p'^k + c p'^l) (\zz^l)^T$. Hence,
$\N( \zz^l (\zz^l)^T)$ = $\N((\zz^l)^T) \subseteq $ $\N( (p'^k+c p'^l) (\zz^l)^T)$.
Therefore, it follows from Lemma \ref{lembasic} that there exists $t \neq  0$ such that
$2B' - t E^{kl} \succeq \bz$, i.e., $t \in \Tl$.

To prove the other direction, assume that $\zz^k$ and $\zz^l$ are not parallel and assume, to the contrary, that
there exists $t \neq 0$ such that $2B' - t E^{kl} \succeq \bz$. Thus
$\N( \zz^k (\zz^l)^T+\zz^l (\zz^k)^T) \subseteq $ $\N( p'^k (\zz^l)^T + p'^l (\zz^k)^T)$.
Next, we consider two cases:

(i)  $\zz^k \neq \bz$ and $\zz^l = \bz$. In this case,
$\zz^k (\zz^l)^T + \zz^l (\zz^k)^T = \bz$ and  $p'^k (\zz^l)^T + p'^l (\zz^k)^T) = p'^l (\zz^k)^T \neq \bz$ since,
by Proposition \ref{prop2}, $p'^l \neq \bz$. Hence, we have a contradiction since $\N(\bz) \not \subseteq \N(p'^k (\zz^l)^T)$.

(ii) both $\zz^k$ and $\zz^l$ are nonzero. Also, in this case we have a contradiction since Proposition \ref{prop1} implies
that $\zz^k (\zz^l)^T+\zz^l (\zz^k)^T$ is indefinite. Consequently, $\Tl=\{0\}$.
\epr

\begin{exa}
Consider the unit spherical EDM
\[
D=\left[ \begin{array}{cccc} 0 & 2 & 4 & 2 \\ 2 & 0 & 2 & 4 \\
                                                              4 & 2 & 0 & 2 \\ 2 & 4 & 2 & 0 \end{array} \right]
\mbox{ and thus }
B^{\dag}= \frac{1}{4} B = \frac{1}{4}  \left[ \begin{array}{rrrr} 1 & 0 & -1 & 0 \\ 0 & 1 & 0 & -1 \\
                                                              -1 & 0 & 1 & 0 \\ 0 & -1 & 0 & 1 \end{array} \right],
\]
$w_1=w_2=w_3=w_4=1/8$ and $z^1=z^3=1$, $z^2=z^4=-1$.
Note that $D$ is a regular EDM of embedding dimension $2$.
Consider the entry $d_{12}$. Clearly, $d_{12}$ is yielding since $z^1=-z^2$.
However, $\zz^1$ is not parallel to $\zz^2$. Moreover,
$||s^1 + s^2||^2 = B^{\dag}_{11}+ B^{\dag}_{22} + 2 B^{\dag}_{12}=1/2$.
Consequently, $[l_{12}, u_{12}]$ = $[0 \; , \;\thc= 8]$ and
$T^{\leq}_{12} = \{0\}$.

On the other hand, consider the entry $d_{13}$. Obviously, $\zz^1 = \zz^3$, i.e., $c=1$, and thus $d_{13}$ is yielding.
Moreover,
 $||s^1 - s^3||^2 = B^{\dag}_{11}+ B^{\dag}_{33} - 2 B^{\dag}_{13}=1$.
Consequently, $[l_{13}, u_{13}]$ = $[\thc= -4 \; , \; 0]$ and $T^{\leq}_{13} \neq \{0\}$.
In fact, as Theorem \ref{thmTln0} below shows, in this case $T^{\leq}_{13} = [l_{13}, u_{13}]$.
\end{exa}

The case where $\Tl \neq \{0\}$ is divided into two subcases depending on whether $\zz^k=\zz^l=\bz$ or whether
$\zz^k= c \zz^l \neq \bz$, where $c$ is a nonzero scalar.
We start by characterizing $\Tl$ in the subcase where $\zz^k=\zz^l=\bz$. As the following theorem shows, in this subcase,
$\Tl$ is equal to the yielding interval $[l_{kl},u_{kl}]$.

\begin{thm} \label{thmTl0int}
Let $D$ be an $n \times n$ unit spherical EDM  of embedding dimension $r$ and let $(\zz^i)^T$ be the $i$th row of $\ZZ$  defined in (\ref{defZZ}).
Assume that  $\zz^k = \zz^l = \bz$, i.e., $w_k=w_l=0$ if $r=n-1$ and $w_k=w_l=0$ and $z^k=z^l=\bz$ if $r \leq n-2$.
Then
 \[
\Tl= \left[ \thl \;, \;\thu \right],
\]
where $\thl$ and $\thu$ are as defined in (\ref{defthul}).
\end{thm}

\bpr
Lemma \ref{lembasic} implies that $t \in \Tl$ if and only if  $2 (P'^TP')^2 - t (p'^k (p'^l)^T + p'^l (p'^k)^T) \succeq \bz$ if and only if
$ 2 I - t(s'^k (s'^l)^T + s'^l (s'^k)^T) \succeq \bz$. The result follows from Propositions \ref{prop2} and \ref{prop1}
and Corollary \ref{corBiBi'w0}.
\epr

Next, we characterize $\Tl$ in the subcase where $\zz^k= c \zz^l \neq \bz$, where $c$ is a nonzero scalar.
As the next theorem shows, in this subcase, $\Tl$ may or may not be equal to the yielding interval $[l_{kl},u_{kl}]$.
We will elaborate on this point after we have proved the following theorem.

\begin{thm} \label{thmTln0}
Let $D$ be an $n \times n$ unit spherical EDM and embedding dimension $r$ and let $(\zz^i)^T$ be the $i$th row of $\ZZ$ defined in (\ref{defZZ}).
Assume that  $\zz^k = c \zz^l \neq \bz$, where $c$ is a nonzero scalar.
Then
 \[
\Tl = \left\{ \begin{array}{ll} \left[  \thc \; , \; 0  \right] & \mbox{ if } c > 0, \\
                                        \left[ 0 \; , \; \thc  \right]  & \mbox{ if } c < 0, \end{array} \right.
\]
where $\thc$ is defined in (\ref{defthc}).
\end{thm}

\bpr
Lemma \ref{lembasic} implies that  $t \in \Tl$ if and only if
\beq \label{eq2}
 \left[ \begin{array}{cc} 2(P'^TP')^2 - t \; (p'^k (p'^l)^T + p'^l (p'^k)^T) &
                   -t \;( p'^k + c p'^l) (\zz^l)^T) \\ -t \; \zz^l (p'^k+c p'^l)^T &
                   - t \; 2 c \zz^l (\zz^l)^T \end{array} \right] \succeq \bz.
\eeq
Let $M$ be a matrix such that $Q'=[ \frac{\zz^l}{||\zz^l||} \; \;  M]$ is an $(n-r) \times (n-r)$
orthogonal matrix. Hence, the $n \times n$ matrix
$Q= \left[ \begin{array}{cc} I_r & \bz \\ \bz & Q' \end{array} \right]$ is obviously orthogonal.
By multiplying the LHS of Equation (\ref{eq2}) from the left with $Q^T$ and from the right with $Q$, we get that
 $t \in \Tl$  if and only if
\beq \label{eq3}
 \left[ \begin{array}{cc} 2(P'^TP')^2 - t \; (p'^k (p'^l)^T + p'^l (p'^k)^T) &
                   -t \;( p'^k + c p'^l) \; ||\zz^l||) \\ -t \; ||\zz^l|| \; (p'^k+c p'^l)^T) &
                   - t \; 2 c \; ||\zz^l||^2 \end{array} \right] \succeq \bz.
\eeq
Using Schur complement, we have that Equation (\ref{eq3}) holds iff
\beq \label{eqnadd}
tc \leq 0 \mbox{ and } 2(P'^TP')^2 + \frac{t}{2c} (p'^k - c p'^l)(p'^k-cp'^l)^T \succeq \bz,
\eeq
which is equivalent to
\[
tc \leq 0 \mbox{ and } 2I_r  + \frac{t}{2c} (s'^k - c s'^l)(s'^k-cs'^l)^T \succeq \bz,
\]
which, in turn, is equivalent to
\[
tc \leq 0 \mbox{ and } 2 + \frac{t}{2c} || s'^k-cs'^l||^2 \geq 0.
\]
The result follows from Proposition \ref{prop2} and from
 Corollary \ref{corBiBi'c}.
\epr

As we mentioned earlier, in Theorem  \ref{thmTln0}, $\Tl$ is equal to the yielding interval $[ l_{kl}, u_{kl}]$
if $r \leq n-2$ and the Gale transform $z^k \neq \bz$. Equivalently, $\Tl$ is a proper subset of $[ l_{kl}, u_{kl}]$ if
$w_k=cw_l \neq 0$ and either $r=n-1$ or $z^k=z^l=\bz$.

\begin{exa} \label{exaw=0}
Consider the unit spherical EDM $D=\left[ \begin{array}{cccc} 0 & 4 & 2 & 2 \\ 4 & 0 & 2 & 2 \\
                                                              2 & 2 & 0 & 2 \\ 2 & 2 & 2 & 0 \end{array} \right]$
of embedding dimension $3$. Then
\[
B= \frac{1}{8} \left[ \begin{array}{rrrr} 9 & -7 & -1 & -1 \\ -7 & 9 & -1 & -1 \\
                                                              -1 & -1 & 5 & -3 \\ -1 & -1 & -3 & 5 \end{array} \right]
\mbox{ and thus }
B^{\dag}= \frac{1}{4} \left[ \begin{array}{rrrr} 3 & 1 & -2 & -2 \\ 1 & 3 & -2 & -2 \\
                                                              -2 & -2 & 4 & 0 \\ -2 & -2 & 0 & 4 \end{array} \right].
\]
Moreover, $w_1=w_2=1/4$ and $w_3=w_4=0$.
Consider the yielding entry $d_{12}$. In this case $c=1$ and thus
\[ ||s^1-cs^2||^2 = B^{\dag}_{11}+ c^2 B^{\dag}_{22}-2cB^{\dag}_{12}=1,
\]
\[
(s^1)^T s^2 - ||s^1||\; ||s^2|| =B^{\dag}_{12} - (B^{\dag}_{11} B^{\dag}_{22})^{1/2} =  -1/2, \mbox{ and }
(s^1)^T s^2 + ||s^1||\; ||s^2|| =1.
\]
Therefore,
$[l_{12}, u_{12}]$ = $[\thl=-4 \; , \;\thu= 2]$, while
$T^{\leq}_{12} = [ \thc= -4 \;,\;0]$. Note that $\thc=\thu$ and
$B^{\dag}_{11} =c^2 B^{\dag}_{22}$.

On the other hand, consider the yielding entry $d_{34}$. Thus,
$(s^3)^T s^4 - ||s^3||\; ||s^4|| = - 1$ and
$(s^3)^T s^4 + ||s^3||\; ||s^4|| = 1$.
Therefore, $[l_{34}, u_{34}]$ =  $T^{\leq}_{34}=[\thl=-2 \; , \;\thu= 2]$.
\end{exa}

\section{Characterizing $\Te$}

Let $d_{kl}$ be a yielding entry of a unit spherical EDM $D$.
Recall that
\[
\Te=\{t \in \Tl  : D+tE^{kl} \mbox{ is a unit spherical EDM} \}.
\]
We saw in the previous section that $\Tl$ is a convex subset of the yielding interval $[l_{kl},u_{kl}]$
of $d_{kl}$. As will be shown in this section, $\Te$ may or may not be equal to $\Tl$. Moreover, in case $\Te \neq \Tl$, $\Te$  may
or may not  be  a convex subset of $\Tl$. We start first with the case where $\Te = \Tl$.

\subsection{The Case Where $\Te=\Tl$}

The equality between $\Te$ and $\Tl$ is proved in the following theorem by establishing a lower bound on
the radius of $D+t E^{kl}$ for $t \in \Tl$. This is achieved by using the duality theory of semidefinite programming (SDP).

\begin{thm} \label{thmrho=1}
Let $D$ be an $n \times n$ unit spherical EDM of embedding dimension $r$ and let $w= D^{\dag} e$.
Let $(\zz^i)^T$ denote the $i$th row of $\ZZ$ defined in (\ref{defZZ}) and assume that $\zz^k=c\zz^l$ for some nonzero scalar $c$.
Further, let $\rho_{kl}(t)$  denote the radius of $D+tE^{kl}$ for
$ t \in \Tl $.  If
\begin{enumerate}[(i)]
\item $w_k=w_l=0$ or
\item $w_k \neq 0$, $r \leq n-2$ and $z^k \neq \bz$,
\end{enumerate}
then
\[
  \rho_{kl}(t) = 1 \mbox{ for all }  t \in \Tl,
\]
and thus $\Te=\Tl$.
\end{thm}

\bpr
Consider the following pair of dual SDP problems
\[
\begin{array}{lll}
\mbox{(P)} &  \min & \lambda                                \\
    & \mbox{subject to }& 2 \lambda E - t E^{kl} \succeq D.
\end{array}
\]
and
\[
\begin{array}{llll}
\mbox{(D)} & \max & \tra(DY) &  \\
    & \mbox{subject to} & 2 \; \tra(EY) & = 1 \\
   &     & Y_{kl}  & = 0 \\
     &    & Y \succeq \bz. &
\end{array}
\]
Since $B'=E - D/2$ where $B'w=\bz$, the dual problem (D) is equivalent to
\[
\begin{array}{lrll}
\mbox{(D)} & 1- \min & 2 \; \tra(B'Y) &  \\
    & \mbox{subject to} & 2 \; \tra(EY) & = 1 \\
   &     & Y_{kl}  & = 0 \\
     &    & Y \succeq \bz. &
\end{array}
\]
Note that $\tra(B'Y) \geq 0$ since both $B'$ and $Y$ are positive semidefinite.
It is easy to see that the Slater's condition  \cite{wsv00} holds for the dual problem (D).
Let $\lambda^*$ and $Y^*$ denote, respectively, the optimal solutions of problems (P) and (D).
Then, by SDP  strong duality, we have
\[
1 \geq \lambda^* = 1 - 2 \; \tra(B'Y^*).
\]
Furthermore,
Theorem \ref{thmsEDM} implies that
the minimum of $\rho^2_{kl}(t) = \lambda^*$.
First, assume that  $w_k=w_l=0$ and let $Y=2 ww^T$. Then $2 e^TYe= 4 (e^Tw)^2 = 1$ and $Y_{kl}=0$ and $Y \succeq \bz$.
Therefore, $Y$ is an optimal solution of (D) since $\tra(B'Y)=0$.
Consequently, $\lambda^* =1$ and thus $\rho^2_{kl}(t) =1$ for all $t \in \Tl$.

Now assume that $w_k \neq 0$, $r \leq n-2$ and $z^k \neq \bz$.
Let $\zeta$ be a column of $Z$ such that $\zeta_k \neq 0$. Such $\zeta$ exists since
$z^k \neq \bz$.  Let $y = w - (w_k / \zeta_k) \zeta$ and let $Y=2 yy^T$. Then $Y \succeq \bz$ and $Y_{kl}=0$ since $y_k=0$.
Moreover, $2 e^T Y e = 4 (e^Tw)^2 =  1$. Therefore, $Y$ is an optimal solution of (D) since $\tra(B'Y)=0$.
Consequently, $\lambda^* =1$ and thus $\rho^2_{kl}(t) =1$ for all $t \in \Tl$.
\epr

As the following subsection shows, the analysis of the case where $\Te \neq \Tl$ is much more involved.

\subsection{The Case Where $\Te \neq \Tl$}

If $\Te \neq \Tl$, then
the following two facts are immediate consequences of Theorem \ref{thmrho=1}:
First, $w_k=c w_l \neq 0$  and either
$r=n-1$ or $z^k=z^l=\bz$. Second,
\[
\min \{ \rho_{kl}(t): t \in \Tl \} < 1,
\]
where $\rho_{kl}(t)$ is the radius of $D+t E^{kl}$.
As a result, we need to find an explicit expression for $\rho_{kl}(t)$.
To this end, let $B=-JDJ/2$ and let $w=D^{\dag}e$. Let us define the following quantities:

\beq \label{defab}
\begin{array}{lll}
\alpha_1 & =  & 2 D^{\dag}_{kl}, \\
 \alpha_2 & = & (D^{\dag}_{kl})^2 - D^{\dag}_{kk} D^{\dag}_{ll}, \\
\beta_1 & = & \alpha_1 - 4 c w_l^2 =- B^{\dag}_{kl}, \\
 \beta_2 & = & \alpha_2 +2 w_l^2( D^{\dag}_{kk}+ c^2 D^{\dag}_{ll} - 2 c D^{\dag}_{kl}) = \frac{1}{4} ((B^{\dag}_{kl})^2-B^{\dag}_{kk} B^{\dag}_{ll});
\end{array}
\eeq
and the following two polynomials:

\beq \label{defft}
f(t)=1 + \alpha_1 t + \alpha_2 t^2 \mbox{ and } g(t)=1+\beta_1 t + \beta_2 t^2.
\eeq

Two remarks are in order here. First, whereas $\alpha_2$ can be zero or nonzero, $\beta_2 < 0$.
To see why this is the case, note that $\beta_2 \leq 0$ since $B^{\dag} \succeq \bz$. Now assume to the
contrary that $\beta_2=0$, then by Cauchy-Schwarz inequality, $p^k$ is parallel to $p^l$.
But this contradicts part 2 of Proposition 2.1 of \cite{alf18}. Second, $g(t)$ can be factorized as
\beq
g(t) = \beta_2(t - \thl) ( t - \thu),
\eeq
where $\thl$ and $\thu$ are as defined in (\ref{defthul}). The following two technical lemmas are critical for the results of this
subsection. Their proofs can be established by straightforward calculations.

\begin{lem} \label{lemfthc}
Let $f(t)$ and $g(t)$ be as defined in (\ref{defft}) and let $\thl$, $\thu$  and $\thc$ be as defined in (\ref{defthul}) and (\ref{defthc}).
Assume that $w_k=c w_l \neq 0$. Then
\[
f(\thl) = \frac{4 w_l^2 \; ( ||s^k|| - c \; ||s^l||)^2 } { ((s^k)^Ts^l-||s^k|| \; ||s^l||)^2 },
\]
\[
f(\thu) = \frac{4 w_l^2 \; ( ||s^k|| + c \; ||s^l||)^2 } { ((s^k)^Ts^l + ||s^k|| \; ||s^l||)^2 },
\]
and
\[
f(\thc) = g(\thc)= \frac{( ||s^k||^2 - c^2 \; ||s^l||^2)^2 } { ||s^k - c s^l||^4 }.
\]
\end{lem}

\begin{lem} \label{lemthc-thu}
Let $\thl$, $\thu$  and $\thc$ be as defined in (\ref{defthul}) and (\ref{defthc}).
Then
\[
\thc-\thl =  \frac{-2 \; ( ||s^k|| - c \; ||s^l||)^2  } {||s^k-cs^l||^2 \; ((s^k)^T s^l - ||s^k|| \; ||s^l||) } \geq 0,
\]
and
\[
\thu-\thc =  \frac{2 \; ( ||s^k|| + c \; ||s^l||)^2  } {||s^k-cs^l||^2 \; ((s^k)^T s^l + ||s^k|| \; ||s^l||) } \geq 0.
\]
\end{lem}

The following theorem, which is the main result of this subsection, provides
a closed-form expression for $\rho_{kl}(t)$.

\begin{thm}  \label{thmrhot}
Let $D$ be an $n \times n$ unit spherical EDM of embedding dimension $r$
and let $B=-JDJ/2$. Let $Z$ be a Gale matrix of $D$ and let $w= D^{\dag} e$.
Further, let $\rho_{kl}(t)$ denote the radius of  $(D + t E^{kl})$, $t \in \Tl$.
Assume that $w_k= c w_l \neq 0$ and either $r=n-1$ or $z^k = z^l = \bz $.
Then
\[
\rho_{kl}^2(t) = \frac{f(t)}{ g(t)}
=  \frac{f(t)}{\beta_2(t - \thl) ( t - \thu)}
\]
for all $t \in \Tl$, where $f(t)$ and $g(t)$ are defined in (\ref{defft}).
\end{thm}

\bpr
Consider first the case where $r=n-1$. Theorem \ref{thmsEDM} implies, in this case, that $D$ is nonsingular.
Let $(D + t E^{kl}) w(t) = e$. Then
$e^T w(t) = e^T ( D + t E^{kl})^{-1} e$.
Making use of the Sherman-Morrison-Woodbury formula twice, we get that
\begin{eqnarray*}
(D+tE^{kl})^{-1} & = &  D^{-1} - \frac{t}{f(t)} \left[  (1 + t D^{-1}_{kl})  (D^{-1}E^{kl} D^{-1} ) \right]\\
                 & & + \frac{t^2}{f(t)} \left[  (D^{-1}_{kk} \; D^{-1}e^l (e^l)^T  D^{-1} + D^{-1}_{ll} \; D^{-1}e^k (e^k)^T  D^{-1} ) \right].
\end{eqnarray*}
Using Theorem \ref{thmsEDM}, the result can be established, in this case, by a straightforward calculation.

Now consider the case where $r \leq n-2$ and let $Z$ be a Gale matrix of $D$. Thus Theorem \ref{thmsEDM} implies that the columns of
$Z$ form a basis of $\N(D)$.
Let $(D + t E^{kl}) w(t) = e$. Then by multiplying this equation from the left with $D^{\dag}$ and using the fact that
$D^{\dag} D = I - Z (Z^TZ)^{-1} Z^T$ we have that
\beq \label{defwt}
w(t) - Z (Z^TZ)^{-1} Z^T w(t) + t D^{\dag} E^{kl} w(t) = w,
\eeq
and hence
\beq \label{eqwt1}
e^T w(t) = e^T w - t w^T E^{kl} w(t).
\eeq
By multiplying (\ref{defwt}) from the left by $E^{kl}$ and noting that by assumption $E^{kl} Z = \bz$, we get that
\beq  \label{eqwt2}
E^{kl} w(t) = (I + t E^{kl} D^{\dag})^{-1} E^{kl}w.
\eeq
Now substituting (\ref{eqwt2}) into (\ref{eqwt1}) we get
\[
e^T w(t) = e^T w - t  w^T (I + t E^{kl} D^{\dag})^{-1} E^{kl} w.
\]
Again the result follows, in this case, by applying the Sherman-Morrison-Woodbury formula twice on  $(I + t E^{kl} D^{\dag})^{-1}$
and by using Theorem \ref{thmsEDM}.
 \epr

It is easy to deduce from Theorem \ref{thmrhot} that if $\rho^2_{kl}(t) = 1$, then $t=0$ or $t=\thc$. Obviously,  $\rho^2_{kl}(0) = 1$.
However, $\rho^2_{kl}(\thc)$ may or may not be equal to 1. The following theorem characterizes the case where
$\rho^2_{kl}(\thc) < 1$.

\begin{thm}\label{thmTe=0}
Let $D$ be an $n \times n$ unit spherical EDM  of embedding dimension $r$
and let $B=-JDJ/2$. Let $Z$ be a Gale matrix of $D$ and let $w= D^{\dag} e$.
Assume that $w_k= c w_l \neq 0$ and either $r=n-1$ or $z^k = z^l = \bz$.
Let $s^i$ be as defined in (\ref{defSS'}) and assume, further, that
$||s^k||^2 = c^2 ||s^l||^2$.  Then
\[
\rho_{kl}(t) = 1 \mbox{ if and only if } t=0,
\]
where $\rho_{kl}(t)$ is the radius of  $(D + t E^{kl})$.
\end{thm}

\bpr
Assume that $c > 0$.
Then it follows from  Lemmas \ref{lemfthc} and \ref{lemthc-thu} that $\thl=\thc$ and $f(\thl)=0$.
Hence, by L'Hospital's rule,
\[
\rho_{kl}^2(\thc) = \frac{f'(\thc)}{g'(\thc)}
  = -2 \frac{D^{\dag}_{ll}}{B^{\dag}_{ll}}=1 - \frac{4 w_l^2}{B^{\dag}_{ll}} < 1.
\]

The case where $c < 0$ implies that $\thu=\thc$ and $f(\thu)=0$ and the proof
is similar to the previous case.
\epr

\begin{thm}  \label{thmTe=0thc}
Let $D$ be an $n \times n$ unit spherical EDM of embedding dimension $r$
and let $B=-JDJ/2$. Let $Z$ be a Gale matrix of $D$ and let $w= D^{\dag} e$.
Assume that $w_k = c w_l \neq 0$ and either $r=n-1$ or $z^k = z^l = \bz $.
Let $s^i$ be as defined in (\ref{defSS'}) and  assume that $||s^k||^2 \neq c^2 ||s^l||^2$.
Then
\[
\rho_{kl}(t) = 1 \mbox{ if and only if } t=0 \mbox{ or } t=\thc,
\]
where $\rho_{kl}(t)$ is the radius of  $(D + t E^{kl})$.
\end{thm}

\bpr This follows from Lemma \ref{lemfthc} and Lemma \ref{lemthc-thu} since
$g(\thc)=f(\thc)$ and $ f(\thc) \neq 0$.
Hence, $\rho_{kl}(\thc) = 1$.
\epr

As a result, we have the following characterization of $\Te$ when it is not equal to $\Tl$.

\begin{cor}
Let $D$ be an $n \times n$ unit spherical EDM  of embedding dimension $r$
and let $B=-JDJ/2$. Let $Z$ be a Gale matrix of $D$ and let $w= D^{\dag} e$.
Assume that $w_k= c w_l \neq 0$ and either $r=n-1$ or $z^k = z^l = \bz$.
Let $s^i$ be as defined in (\ref{defSS'}).
Then
\[
\Te = \left\{ \begin{array}{ll} \{0\} & \mbox{ if } ||s^k||^2 = c^2 ||s^l||^2, \\
                               \{0,\thc\} & \mbox{ otherwise. } \end{array} \right.
\]
\end{cor}

\begin{exa}\label{exalast}
Consider the EDM $D = \left[ \begin{array}{ccc} 0 & 1 & 3 \\ 1 & 0 & 1 \\ 3 & 1 & 0 \end{array} \right]$.
Then $D$ is unit spherical of embedding dimension $2$ and
$w = \frac{1}{2}[ 1 \; -1 \; 1]^T$. Moreover,
a configuration matrix of $D$ is
\[ P=\frac{1}{6} \left[ \begin{array}{cr} -3 \sqrt{3} & 1 \\ 0 & -2 \\ 3 \sqrt{3} & 1 \end{array} \right]
\mbox{ and thus }
S= \left[ \begin{array}{cr} -1/ \sqrt{3} & 1 \\ 0 & -2 \\ 1/ \sqrt{3} & 1 \end{array} \right].
\]
Consider the entry $d_{12}$.
The yielding interval of $d_{12}$ is $[\thl= 3 - 2 \sqrt{3}, \thu= 3 + 2 \sqrt{3}]$ and
$T^{\leq}_{12} = [0 , \thc= 3]$.
Moreover, $||s^1|| \neq ||s^2||$. Therefore,
\[
\rho_{12}^2(t) = \frac{3+ 3t}{3 + 6 t - t^2},
\]
and hence $\rho^2_{12}(0) = \rho^2_{12}(3) = 1$. Consequently, $T^{=}_{12}=\{0,3\}$.

On the other hand, consider the entry $d_{13}$.
The yielding interval of $d_{13}$ is  $[\thl=-3, \thu=1]$ while
$T^{\leq}_{13} = [\thc= -3 , 0]$. Note that in this case
 $||s^1|| = ||s^2||$ and thus $\thc=\thl$.
Furthermore.
\[
\rho_{13}^2(t) = \frac{1}{(1-t)}.
\]
Consequently,  $T^{=}_{13}=\{0\}$. It is worth pointing out that for
$0 < t < \thu$, $D+t E^{13}$ is a spherical EDM of radius $\rho > 1$; and
$D+ \thu E^{13}$ is a nonspherical EDM.
\end{exa}

\section{Alternative Approach}

The results in Sections 3 and 4 can be alternatively derived by using the Cayley-Menger matrix.
Given an $n \times n$ EDM $D$, the $(n+1) \times (n+1)$ matrix
\beq
\DD = \left[ \begin{array}{cc} 0 & e^T \\ e & D \end{array} \right]
\eeq
is called the Cayley-Menger matrix of $D$.
Let $\ee$ denote the vector of all 1's in $\Rs^{n+1}$.
Let $U= \left[ \begin{array}{c} - e^T \\ I_n \end{array} \right]$. Then obvisously,
$U(U^TU)^{-1}U^T$ is a projection matrix on $\ee^{\perp}$. Consequently, the Cayley-Menger matrix
$\DD$ is an EDM iff ($- U^T \DD U \succeq \bz$), i.e., iff
\beq \label{eqDDEDM}
2E - D \succeq \bz.
\eeq

The following theorems establish some relations between $D$ and $\DD$.

\begin{thm} \label{lemDDrho}
Let $\DD$ be the Cayley-Menger matrix of $D$. Then $\DD$ is an EDM
if and only if $D$ is a spherical EDM of radius $\rho \leq 1$, in which case,
$\rho^2 = 1- e^T \ww/2$ where  $\DD \ww = \ee$.
\end{thm}

\bpr
It follows from (\ref{eqDDEDM}) and  Theorem \ref{thmsEDM} that
$\DD$ is an EDM if and only if $2E - D \succeq \bz$ if and only if
 $D$ is an EDM of radius $\rho \leq 1$.

Now assume that $\DD$ is an EDM and let $\DD \ww = \ee$ and $Dw=e$. Then it is easy to show that
$\ww= \left[ \begin{array}{c} 1- 2 \rho^2 \\ 2 \rho^2 w \end{array} \right]$.
Consequently, $\ee^T \ww = 2 - 2 \rho^2$ and the result follows.
\epr

The following theorem is an immediate consequence of Theorem \ref{lemDDrho}.
\begin{thm}
$D$ is a unit spherical EDM
if and only if its Cayley-Menger matrix $\DD$  is a nonspherical EDM.
\end{thm}

\bpr Assume that $\DD$ is a nonspherical EDM and let $\DD \ww = \ee$. Then by Theorem \ref{thmnsEDM}, $\ee^T \ww = 0$.
Hence, by Theorem \ref{lemDDrho}, $D$ is a unit spherical EDM.

To prove the other direction, assume that $D$ is a unit spherical EDM and let $Dw=e$.
Then $\DD$ is an EDM. Moreover, let $\DD \ww = \ee$. Then
$\ww= \left[ \begin{array}{c} -1 \\ 2 w \end{array} \right]$. Hence,
$\ee^T \ww = 0$ since $e^Tw = 1/2$.
\epr

The following lemma establishes the equality of the embedding dimensions
of $D$ and $\DD$ when $D$ is a unit spherical EDM.

\begin{lem}  \label{lemr+2}
Let $D$ be a unit spherical EDM and of embedding dimension $r$.  Let
$\DD$ be the Cayley-Menger matrix of $D$. Then
the embedding dimension of $\DD$ is $r$.
\end{lem}

\bpr
Let $Dw=e$ and let $Q= \left[ \begin{array}{ccc} 1  & 0 & \bz \\ \bz & w & V \end{array} \right]$.
Then $Q$ is nonsingular since $e^Tw=1/2 > 0$. Moreover,
\[
Q^T \DD Q = \left[ \begin{array}{ccc} 0 & e^Tw & \bz \\ e^Tw & e^Tw & \bz \\ \bz & \bz & V^TDV \end{array} \right].
\]
Consequently, $\rank(\DD) = r+2$ and thus, by Theorem \ref{thmsEDM},
 the emedding dimension of $\DD$ is $r$ since $\DD$ is nonspherical.
\epr

Next, we establish the relations between configuration and Gale matrices of $D$ and $\DD$.
Let $D$ be a unit spherical EDM of center $a$ and let $P$ be a configuration of $D$ such that $P^Te=\bz$. Then
it is not difficult to show that
$\PP =  \left[ \begin{array}{l} a^T \\ P \end{array} \right]$ is a configuration matrix of $\DD$.
Moreover,

\begin{lem}  \label{lemZDZZDD}
Let $D$ be an $n \times n$ unit spherical EDM of embedding dimension $r$.
Let $Dw=e$. If $r=n-1$, then
$\ZZ =  \left[ \begin{array}{c} -1/2 \\ w \end{array} \right]$ is a Gale matrix of $\DD$.
Otherwise, i.e., if $r \leq n-2$, then
$\ZZ =  \left[ \begin{array}{cc} -1/2 & \bz \\ w & Z \end{array} \right]$ is a Gale matrix of $\DD$,
where $Z$ is a Gale matrix of $D$.
\end{lem}

\bpr
Since $ e^T w =1/2$,
it suffices to show that $2 P^T w = a$. To this end, we have
\[
e = Dw = (\diag(B))^T w \; e + e^Tw \; \diag(B) - 2 B w.
\]
Thus, $2 P^TBw = P^T \diag(B)/2$. But, by Theorem \ref{thmsEDM},
 $2 Pa = J \diag(B)$ and thus $2 P^T P a = P^T \diag(B)$. Hence,
 $2 P^TP P^Tw = P^TP  a$ and thus the result follows.
\epr

As a result,
using Lemma \ref{lemZDZZDD} and keeping in mind that $2 \leq k < l \leq n+1$,
Theorem \ref{thmTl0} follows by applying Theorem \ref{thmfrom2} to $\DD$. Also,
using Lemma \ref{lemZDZZDD},
Theorems \ref{thmTl0int} and \ref{thmTln0} follow by applying Theorem \ref{thmfrom2int} to $\DD$.

The following theorem is crucial to the rederivation of  the results concerning $\Te$.

\begin{thm} \label{thmDDew}
Let $D$ be a unit spherical EDM of embedding dimension $r$.
Let $\DD(t)$ be the Cayeley-Menger matrix of $D+t E^{kl}$. Let $w= D^{\dag} e$ and let
$\DD(t) \ww(t) = \ee$.
Assume that $w_k = c w_l \neq 0$ and either $r=n-1$  or $z^k=z^l=\bz $.
Then
\[
\ee^T \ww(t) =\frac{8 w_l^2 c \, t}{\thc \beta_2}  \times
\left\{ \begin{array}{ll} {\displaystyle \frac{(t - \thc) }{(t-\thl)(t-\thu)} } & \mbox{ if } \thc \neq \thl \mbox{ and } \thc \neq \thu, \\
                          {\displaystyle          \frac{1 }{(t-\thu)}  }  & \mbox{ if } \thc = \thl, \\
                            {\displaystyle           \frac{1 }{(t-\thl)} } & \mbox{ if } \thc = \thu. \end{array} \right.
\]
\end{thm}

\bpr
Let $\ww(t) = \left[ \begin{array}{c} \sigma(t) \\ \xi(t) \end{array} \right]$.
Then
\beq  \label{eqz=0diff0}
\left[ \begin{array}{cc} 0 & e^T \\ e & D \end{array} \right]
\left[ \begin{array}{c} \sigma(t) \\ \xi(t) \end{array} \right]
+ t
\left[ \begin{array}{cc} 0 & \bz \\ \bz & E^{kl} \end{array} \right]
\left[ \begin{array}{c} \sigma(t) \\ \xi(t) \end{array} \right] =
\left[ \begin{array}{c} 1 \\ e \end{array} \right].
\eeq
Assume first that $r=n-1$, i.e., $D$ is nonsingular. Then, by Lemma \ref{lemr+2},
$\DD$, the Cayley-Menger matrix of $D$, is also nonsingular.
Using Lemma \ref{lemDdagBdag}, we
multiply (\ref{eqz=0diff0}) from the left by $\DD^{-1} =
\left[ \begin{array}{cc} -2 & 2 w^T \\ 2 w & -  B^{\dag}/2 \end{array} \right]$ to get
that
\beq \label{eqz=0diff20}
\left[ \begin{array}{c} \sigma(t) \\ \xi(t) \end{array} \right] +
 \left[ \begin{array}{cc} 0 & 2t w^T E^{kl} \\ \bz & - t B^{\dag}E^{kl}/2 \end{array} \right]
\left[ \begin{array}{c} \sigma(t) \\ \xi(t) \end{array} \right] =
\left[ \begin{array}{c} -1 \\ 2w \end{array} \right].
\eeq

Therefore,
\[
\ee^T \ww(t) = - 2 t w^T E^{kl} \xi(t).
\]
Now multiplying  (\ref{eqz=0diff20}) from the left by
$\left[ \begin{array}{cc} 0 & \bz \\ \bz & E^{kl} \end{array} \right]$, we get that
\[
E^{kl} \xi(t) = 2 (I - t E^{kl} B^{\dag}/2)^{-1} E^{kl} w.
\]
Consequently,
\[
\ee^T \ww(t) = - 4 t w^T (I - t E^{kl} B^{\dag}/2)^{-1} E^{kl} w.
\]
But,

\begin{eqnarray*}
(I - t E^{kl} B^{\dag}/2)^{-1}  =   I & + &  \frac{t^2 (B^{\dag}_{kk} e^l (e^l)^T + B^{\dag}_{ll}  e^k (e^k)^T) B^{\dag} }{4 g(t)}  \\
& & + \frac{(2 t - B^{\dag}_{kl}t^2)  E^{kl} B^{\dag}  }{4 g(t)},
\end{eqnarray*}

where $g(t)$ is as defined in (\ref{defft}).
Therefore,
\[
w^T (I - t E^{kl} B^{\dag}/2)^{-1} E^{kl} w =  w_l^2 \frac{ (4c + t \, ||s^k - c s^l ||^2 ) }{2 g(t)}.
\]
Consequently,
\[
\ee^T \ww(t) = - 2 t w_l^2 \frac{( 4c + t \, ||s^k - c s^l ||^2) }{g(t)},
\]
and the result follows in this case.

Now assume that $r \leq n-2$ and let $Z$ be a Gale matrix of $D$.
Using Lemma \ref{lemDdagBdag}, we
multiply (\ref{eqz=0diff0}) from the left by $\DD^{\dag} =
\left[ \begin{array}{cc} -2 & 2 w^T \\ 2 w & -  B^{\dag}/2 \end{array} \right]$
and we use the fact that
\[
\DD^{\dag} \DD = \left[ \begin{array}{cc} 1 & \bz \\ \bz &  D^{\dag} D \end{array} \right] =
\left[ \begin{array}{cc} 1 & \bz \\ \bz &  I - Z(Z^TZ)^{-1}Z^T \end{array} \right]
\]
to get
\[
\left[ \begin{array}{c} \sigma(t) \\ \xi(t) \end{array} \right] +
 \left[ \begin{array}{cc} 0 & 2t w^T E^{kl} \\ \bz & - Z(Z^TZ)^{-1}Z^T - t B^{\dag}E^{kl}/2 \end{array} \right]
\left[ \begin{array}{c} \sigma(t) \\ \xi(t) \end{array} \right] =
\left[ \begin{array}{c} -1 \\ 2w \end{array} \right].
\]
The proof proceeds as in the previous case by using the fact that $E^{kl} Z = \bz$.
\epr

As a result, if  $||s^k||^2 = c^2 ||s^l||^2$, then, by Lemma \ref{lemthc-thu},  either $\thc = \thl$ (if $c >0$) or
$\thc=\thu$ (if $c < 0$ ). Thus,
Theorem \ref{thmDDew} implies that $\ee^T\ww(t) = 0$ if and only if $t=0$. On the other hand, if
$||s^k||^2 \neq c^2 ||s^l||^2$, then,  by Lemma \ref{lemthc-thu}, $\thc \neq \thl$ and $\thc \neq \thu$. Therefore,
Theorem \ref{thmDDew} implies that $\ee^T\ww(t) = 0$ if and only if $t=0$ or $t=\thc$.
In other words,
 Theorem \ref{thmDDew} and Lemma \ref{lemthc-thu} imply Theorems \ref{thmTe=0} and \ref{thmTe=0thc}.

Finally, by combining Lemma \ref{lemDDrho} and Theorem \ref{thmDDew} and under the premise of Theorem \ref{thmDDew}, we have that
\beq  \label{eq2rho}
\rho^2_{kl}(t) = 1- \frac{1}{2} \ee^T \ww(t) = 1 -\frac{4 w_l^2 c t(t-\thc)}{\thc g(t)}.
\eeq
But $g(t) - 4 w_l^2c t(t-\thc)/\thc = f(t)$. Thus Equation (\ref{eq2rho}) is identical to that in Theorem \ref{thmrhot}.

\begin{exa}
Let $D$ be the unit spherical EDM considered in Example \ref{exalast}. Then for entry $d_{12}$, we have
$c=-1$, $w_2= -1/2$, $\thc=3$ and $g(t) = 1 + 2 t - t^2/3$. Thus
\[
\rho^2_{12}(t)= 1 - \frac{-t (t-3)}{3 + 6 t - t^2} =  \frac{3+3t}{3 + 6 t - t^2}
\]
as was obtained in Example \ref{exalast}.
On the other hand, for enrty $d_{13}$, we have
$c=1$, $w_2= 1/2$, $\thc=-3$ and $g(t) = 1 - 2 t/3 - t^2/3$. Thus
\[
\rho^2_{13}(t)= 1 - \frac{t (t+3)}{-3 + 2 t + t^2} =  \frac{1}{1- t}
\]
as was obtained in Example \ref{exalast}
\end{exa}


\begin{thebibliography}{10}

\bibitem{alfm18}
A.~Y. Alfakih.
\newblock {\em Euclidean distance matrices and their applications in rigidity
  theory}.
\newblock Springer, 2018.

\bibitem{alf18}
A.~Y. Alfakih.
\newblock On yielding and jointly yielding entries of {E}uclidean distance
  matrices.
\newblock {\em Linear Algebra Appl.}, 556:144--161, 2018.

\bibitem{aw02}
A.~Y. Alfakih and H.~Wolkowicz.
\newblock Two theorems on {E}uclidean distance matrices and {G}ale transform.
\newblock {\em Linear Algebra Appl.}, 340:149--154, 2002.

\bibitem{bb07}
R.~Balaji and R.~B. Bapat.
\newblock On {E}uclidean distance matrices.
\newblock {\em Linear Algebra Appl.}, 424:108--117, 2007.

\bibitem{cri88}
F.~Critchley.
\newblock On certain linear mappings between inner-product and squared distance
  matrices.
\newblock {\em Linear Algebra Appl.}, 105:91--107, 1988.

\bibitem{gal56}
D.~Gale.
\newblock Neighboring vertices on a convex polyhedron.
\newblock In {\em Linear inequalities and related system}, pages 255--263.
  Princeton University Press, 1956.

\bibitem{gow82}
J.~C. Gower.
\newblock Euclidean distance geometry.
\newblock {\em Math. Sci.}, 7:1--14, 1982.

\bibitem{gow85}
J.~C. Gower.
\newblock Properties of {E}uclidean and non-{E}uclidean distance matrices.
\newblock {\em Linear Algebra Appl.}, 67:81--97, 1985.

\bibitem{gru67}
B.~Gr{\"{u}}nbaum.
\newblock {\em Convex polytopes}.
\newblock John Wiley \& Sons, 1967.

\bibitem{ht93}
T.~L. Hayden and P.~Tarazaga.
\newblock Distance matrices and regular figures.
\newblock {\em Linear Algebra Appl.}, 195:9--16, 1993.

\bibitem{neu81}
A.~Neumaier.
\newblock Distance matrices, dimension and conference graphs.
\newblock {\em Nederl. Akad. Wetensch. Indag. Math.}, 43:385--391, 1981.

\bibitem{sch35}
I.~J. Schoenberg.
\newblock Remarks to {M}aurice {F}r\'{e}chet's article: Sur la d\'{e}finition
  axiomatique d'une classe d'espaces vectoriels distanci\'{e}s applicables
  vectoriellement sur l'espace de {H}ilbert.
\newblock {\em Ann. Math.}, 36:724--732, 1935.

\bibitem{ss97}
G.~P.~H. Styan and G.~E. Subak-Sharpe.
\newblock Inequalities and equalities associated with the {C}ampbell-{Y}oula
  generalized inverse of the indefinite admittance matrix of resistive
  networks.
\newblock {\em Linear Algebra Appl.}, 250:349--370, 1997.

\bibitem{thw96}
P.~Tarazaga, T.~L. Hayden, and J.~Wells.
\newblock Circum-{E}uclidean distance matrices and faces.
\newblock {\em Linear Algebra Appl.}, 232:77--96, 1996.

\bibitem{wsv00}
H.~Wolkowicz, R.~Saigal, and L.~Vandenberghe, editors.
\newblock {\em Handbook of Semidefinite Programming. Theory, Algorithms and
  Applications}.
\newblock Kluwer Academic Publishers, Boston MA, 2000.

\bibitem{yh38}
G.~Young and A.~S. Householder.
\newblock Discussion of a set of points in terms of their mutual distances.
\newblock {\em Psychometrika}, 3:19--22, 1938.

\end{thebibliography}

\end{document}